**UN PROBLEMA DA DISCUTERE**

*un inganno con le percentuali*

> *Il peso netto di alcune casse di cocomeri è* 100 kg; *i cocomeri sono formati al* 99% *da acqua. Dopo un certo tempo, in cui a causa della temperatura estiva parte dell'acqua è evaporata, l'acqua costituisce il* 98% *dei cocomeri. Qual è il peso finale dei cocomeri?*

Il problema, tratto da [2], coinvolge solo concetti elementari e, in effetti, può essere proposto in una Scuola secondaria di I grado. Ma è facile sbagliare, anche per un matematico con una certa esperienza, perché la risposta corretta è decisamente poco intuitiva. Le percentuali riservano sempre qualche sorpresa (si veda anche [1]).

Intuitivamente è spontaneo un ragionamento del tipo seguente. All'inizio l'acqua costituisce il 99% dei cocomeri e quindi corrisponde a un peso di 99 kg; alla fine la percentuale dell'acqua è passata al 98% e quindi avremo 98 kg di acqua. Il resto è rimasto inalterato; in conclusione il peso finale è di 99 kg.

Tuttavia, basta una semplice verifica per rendersi conto che questo risultato è sbagliato. Se avessimo 98 kg di acqua su un totale di 99 kg, allora la percentuale dell'acqua sarebbe di pochissimo inferiore al 99%: infatti, la divisione 98/99 dà come quoziente il numero periodico 0,989898...

Forse, il risultato corretto risulta meno strano se facciamo riferimento alla parte dei cocomeri che non è acqua, credo si chiami *residuo secco* (o *sostanza secca*). Questo all'inizio era l'1 % e alla fine è il 2%, cioè è raddoppiato; quindi ...

Vediamo allora di risolvere il problema non cercando una soluzione intuitiva ma applicando i classici metodi dell'algebra elementare. Chiamiamo $x$ il peso finale. Visto che il residuo secco è di 1 kg, il peso dell'acqua alla fine è di $(x-1)$ kg; abbiamo così l'equazione

$$\frac{x-1}{x} = \frac{98}{100}$$

e si trova senza alcuna difficoltà che la risposta è 50 kg, cioè che il peso dei cocomeri si è dimezzato!

Vediamo un altro metodo risolutivo, sostanzialmente equivalente, leggermente più lungo ma forse più chiaro. Questa volta indichiamo con $x$ la quantità d'acqua evaporata (in kg). Si sa che

$$\frac{\text{peso iniziale acqua}}{\text{peso iniziale cocomeri}} = \frac{99}{100}$$

da cui 100×(peso iniziale acqua) = 99×(peso iniziale cocomeri). Dopo l'evaporazione si ha

$$\frac{\text{peso iniziale acqua} - x}{\text{peso iniziale cocomeri} - x} = \frac{98}{100}$$

Con semplicissimi calcoli, ricordando l'uguaglianza precedente, si trova ancora che $x$ è uguale alla metà del peso iniziale dei cocomeri: dopo l'evaporazione il peso dei cocomeri si è ridotto alla metà.

Cerchiamo di visualizzare la situazione con una figura. Per rendere la figura più "leggibile", suppongo che l'acqua iniziale sia il 90% dei cocomeri e che, dopo l'evaporazione si sia ridotta all'80%. Considero un rettangolo (che corrisponde alla partita di cocomeri) e indico il residuo secco con il colore, lasciando in grigio la parte che corrisponde all'acqua (figura 1).

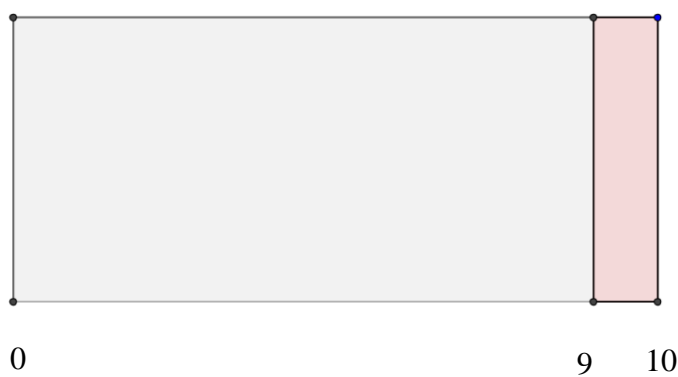

figura 1

La parte grigia è il 90% dell'intero rettangolo. Quanto si deve tagliare a sinistra se vogliamo che, ferma restando la parte colorata, la parte grigia sia l'80% del totale?

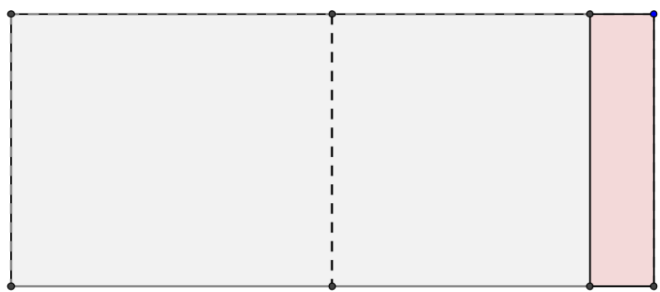

figura 2

Guardando con attenzione la figura 2, è abbastanza chiaro che il taglio a sinistra di una piccola fetta del rettangolo non altera in maniera significativa il rapporto fra la parte colorata e l'intero rettangolo; invece la linea tratteggiata, che dimezza il rettangolo, risponde alla richiesta.

Per capire meglio il legame fra le varie quantità in gioco, indichiamo con $A$ la riduzione dell'acqua in punti percentuali (nel problema proposto $A = 1$) e con $P$ il peso finale in kg dei cocomeri. È chiaro che:

per $A = 0$ (in assenza di evaporazione) si ha $P = 100$,

mentre per $A = 99$ (tutta l'acqua è evaporata) si ha $P = 1$.

Le due righe precedenti potrebbero erroneamente far supporre che $P = 100 - A$ (e nel caso $A = 1$ si troverebbe davvero $P = 99$). Con calcoli analoghi a quelli visti si trova invece

$$P = \frac{100}{1+A}$$

Così, per $A = 1$ si ottiene $P = 50$.

Concludo proponendo due problemi del tutto analoghi.

Nel primo caso, ambientato in un contesto molto diverso, considero volutamente numeri piccoli, che permettano semplici verifiche. Anche in questo caso è probabile che la prima risposta sia sbagliata, ma basta un po' di riflessione per ottenere la risposta corretta.

> *Una classe di 30 studenti è formata da maschi per i 2/3. In seguito a una riorganizzazione delle classi, alcuni maschi vengono spostati in un'altra sezione (e nessuno studente viene aggiunto nella classe considerata); quella classe risulta ora formata da maschi solo per 1/3. Quanti sono, alla fine, gli studenti della classe?*

La seconda variante è più simile al problema iniziale, ma il passaggio è inverso, nel senso che si passa dal 98% al 99%. Forse, in questa forma la risposta corretta è più facilmente accettabile.

> *Una soluzione farmaceutica è formata per il 98% da un eccipiente (acqua), mentre il principio attivo è il 2%. Un flacone è di 500 g. Se si decide di diluire la soluzione, in modo che la percentuale dell'acqua salga al 99%, quanta acqua occorre aggiungere?*

Claudio Bernardi
claudio.bernardi@uniroma1.it

**Riferimenti bibliografici**

[1]   C. Bernardi, *Un problema da discutere – Le percentuali*, Archimede, **LVI** (2004), pag. 93-96

[2]   F. Lazebnik, *Surprises*, Mathematics Magazine, Vol. 87, No. 3 (June 2014), pag. 212-221